\documentclass{kybernetika}

\usepackage{pmat}    

\newtheorem{theorem}{Theorem}[section]
\newtheorem{lemma}[theorem]{Lemma}

\newtheorem{remark}[theorem]{Remark}

\newtheorem{definition}[theorem]{Definition}

\newtheorem{assumption}[theorem]{Assumption}

\begin{document}
\pagestyle{myheadings}

\title{Leader-Following Consensus of Multiple Linear Systems Under Switching Topologies: An Averaging Method}

\author{Wei Ni, Xiaoli Wang and Chun Xiong}

\contact{Wei}{Ni}{School  of Science, Nanchang University, Nanchang 330031,
P. R.  China.}{niw@amss.ac.cn}
\contact{Xiaoli}{Wang}{School of Information Science and Engineering,
Harbin Institute of Technology at Weihai,
  Weihai 264209,  P. R. China.}{xiaoliwang@amss.ac.cn}
\contact{Chun}{Xiong}{School  of Science, Nanchang University, Nanchang 330031,
P. R.  China.}{grigxc@yahoo.com.cn}

\markboth{W. Ni, X. Wang and C. Xiong} {Leader-Following Consensus of Multiple Linear Systems}

\maketitle

\begin{abstract}
The leader-following consensus of multiple linear time invariant (LTI) systems under switching topology is considered.
The leader-following consensus problem consists of designing  for each agent a distributed protocol to make  all agents track a leader vehicle, which has the same LTI dynamics as the agents.
The interaction topology describing the information exchange of these agents is time-varying. An averaging method is
proposed. Unlike the existing  results in the literatures which assume the LTI agents to be  neutrally stable, we relax this condition, only making  assumption that the LTI agents are stablizable and detectable. Observer-based leader-following consensus is also considered.
\end{abstract}



\section{Introduction}
Multi-agent system is a hot topic in a variety of research communities, such as  robotics,
sensor networks, artificial intelligence, automatic control and  biology.
Of particular interest in this field is the consensus problem, since it lays foundations for many consensus-related problem, including formation, flocking  and swarming.  We refer to survey papers \cite{olfati2007,ren2007} and
references therein for details.

Integrator and double integrator models are the simplest abstraction, upon which a large part of results on consensus
of  multi-agent systems have been based (see \cite{ren2005,olfati2004,olfati2007,jadb2003,cheng2008,hong2007}). To deal with more complex models, a number of recent papers are devoted
 to consensus of multiple LTI systems
 \cite{zhang2011,wang2008,ni2010,scardovi2009,seo2009,liu2009,khoo2009,Yoshioka2008,namerikawa2008,wang2009,wang2010,wang2011}.
These results keep most of the concepts provided by earlier developments, and provide new design and analysis technique,
such as LQR approach, low gain approach, $H_{\infty}$ approach, parametrization and geometric approach,  output regulation approach, and  homotopy based approach. However, most of these results
\cite{zhang2011,wang2008,ni2010,seo2009,liu2009,khoo2009,Yoshioka2008,namerikawa2008} mainly focus on fixed interaction topology, rather than time-varying
topology.  How the switches of the interaction topology and agent dynamics jointly affect
the collective behavior of the multi-agent system?  Attempts to understand this issue  had been hampered by the lack of suitable analysis tools.
The results of Scardovi et al. \cite{scardovi2009} and Ni et al. \cite{ni2010} are  mentioned here,
because of their contributions to dealing with switching topology in the setup of high-order agent model. However, when dealing with switching topology,
 \cite{scardovi2009} and \cite{ni2010}  assumed  that  the system of each agent is neutrally stable; thus it has no positive real parts eigenvalues.
This assumption was widely assumed  in the literatures  when the interaction topology is fixed or switching.   Unfortunately,  when the agent is stabilizable and detectable rather than neutrally stable, and when the interaction topology is switching,  there is no
 result reported in the literature to investigate the  consensus of these agents.

To deal with switching  graph topology and to remove the neutral stability condition ,  we provide a  modified averaging approach, which is motivated by \cite{aeyels1999,bellman1985,kosut1987}.
The averaging approach was initially proposed  by Krylov and Bogoliubov in celestial mechanics \cite{krylov1943}, and
was further developed in the work of  \cite{bogoliubov1961,krasnosel1955}; for more details refer to the recent book \cite{sanders2007}. Closely related to the averaging theory is the stability of fast time-varying nonautonomous systems \cite{aeyels1999,kosut1987,bellman1985}, and more specifically the fast switching systems \cite{stilwell2006,teel2011}. The modified  approach in this paper is motivated by the work of Stilwell et al. \cite{stilwell2006}, and also the work of
\cite{aeyels1999,kosut1987,bellman1985}. Although this work borrows the idea from \cite{stilwell2006}, the main difference of this work from \cite{stilwell2006} is  as follows. The synchronization in \cite{stilwell2006} is achieved under under fast switching condition; that is, synchronization is realized under two time scales: a time scale $t$ for the agent dynamics and a time scale for the switching signal parameterized by $t/\varepsilon$ with $\varepsilon$ small enough. In our paper, we further establish that the two time scales can be made to be the same and thus the consensus in our paper is not limited to fast switching case. Furthermore, We present an extended averaging approach for consensus: a sequence of averaged systems (rather a single averaged systems) are the indicators of consensus of the the multi-agent systems. This allows to obtain more relax  conditions for consensus. At last, We give further investigation on how to render these sequence of averaged systems to achieve consensus, and thus ensure the consensus of the original multi-agent systems. This was not investigated in [22]. The result in our paper shows that  if there exists an infinite sequence of uniformly bounded and contiguous time intervals such that during each such interval the interaction graph is jointly connected  and if the dwell time for each subgraph is appropriately small, then consensus can be achieved.

In summary, the contributions of this paper are  as follows:
  \begin{itemize}
\item  Averaging method  is applied to  leader-following consensus of multiple LTI systems.
\item  Results are obtained for a wider class of agent dynamics which is stabilizable and detectable than the existing class of neutrally stable agent dynamics.
\item The agent dynamics and the switching signal considered in this paper have the same time scale, rather than having different time scales considered in \cite{stilwell2006}. Thus the results in our paper are not limited to fast time switching case.
\end{itemize}

The rest of this paper is organized as follows.
Section 2 contains the problem formulation  and some preliminary results.
Section 3 provides the main result of leader-following consensus, and extensions are made in
Section 4 which devotes to   observer-based protocols design and analysis.
Two illustrated examples are presented in Section 5.
Section 6 is a brief conclusion.

\section{Problem Formulation and Preliminaries}
This section presents the multi-agent system model, with each agent being a stabilizable LTI system, which includes integrator or double integrator as its special cases.   The leader-following consensus problem is formulated by use of the graph theory. Some supporting lemmas are also included here.

Consider $N$ agents with the same dynamics
\begin{eqnarray}\label{2.1}
\dot x_i=Ax_i+Bu_i,\quad i=1, 2, \cdots, N,
\end{eqnarray}
where $x_i\in \mathbb{R}^n$ is the agent $i$'s state, and $u_i\in
\mathbb{R}^m$ is agent $i$'s input  through which the interactions
or coupling between agent $i$ and other agents are realized.
 The matrix $B$ is of full column rank.
The state information is transmitted among these agents, and the agents together with the transmission channels form a network.
We use a directed graph $\mathcal {G} =(\mathcal{V}, \mathcal {E})$ to describe the topology of this network,  where $\mathcal{V}=\{1,2,\cdots,N\}$ is the set of nodes
representing $N$ agents and $\mathcal{E} \subset \mathcal{V} \times \mathcal{V}$ is the set of ordered
edges $(i, j)$, meaning that agents $i$ can send information to agent $j$.

The leader, labeled as $i=0$,  has linear dynamics as
\begin{eqnarray}\label{2.2}
\dot x_0=Ax_0,
\end{eqnarray}
where $x_0\in \mathbb{R}^n$ is the state.
Referring agent $i\in \{1, \cdots, N\}$ as follower agent, the leader's dynamics is   obviously independent of follower agents. More specifically, the leader just sends information to some follower agents, without receiving information from them.
The interaction structure of the whole agents $\{0,1,\cdots,N\}$ is described by an extended directed graph
$\bar{\mathcal{G}}=(\bar{\mathcal{V}}, \bar{\mathcal{E}})$, which consists of graph
$\mathcal{G}$, node $0$ and  directed edges from  node $0$ to its information-sending  follower nodes.

\begin{definition} {\bf (Definitions related  to graph)}
Consider a graph $\bar{\mathcal{G}}=(\bar{ \mathcal{V}}, \bar{ \mathcal{E}})$ with
$\bar{ \mathcal{V}}=\{0, 1, \cdots, N\}$.
\begin{itemize}
\item  The set of neighbors of node $i$ relative to subgraph $\mathcal{G}$ of $\bar{\mathcal{G}}$
is denoted by $\mathcal{N}_i=\{j\in \mathcal{V}: (j,i)\in \mathcal{E}, j\neq i\}$.
\item  A directed path is a sequence of  edges $(i_1,i_2),(i_2,i_3),(i_3,i_4),\cdots$ in that graph.
\item  Node $i$ is reachable to node $j$ if there is a directed path from $i$ to $j$.
\item The graph $\bar{\mathcal{G}}$ is called connected if node $0$ is reachable to any other node.
\end{itemize}
\end{definition}


\begin{definition}{\bf (Structure matrices of graph)}\label{def2}
\begin{itemize}
\item For a directed graph $\mathcal{G}$ on nodes
$\{1,\cdots,N\}$, it structure is described by its adjacency matrix $\mathcal {A}\in \mathbb{R}^{N\times N}$,
 whose $ij$-th entry is 1 if
 $(j,i)$ is an edge  of $\mathcal{G}$ and 0 if it is not; or by its Laplacian matrix $\mathcal{L}=-\mathcal {A}+ \Lambda$, where
$\Lambda\in \mathbb{R}^{N\times N}$ is the in-degree matrix of $\mathcal{G}$ which is  diagonal with $i$-th diagonal element be $|\mathcal{N}_i|$, the cardinal of  $\mathcal{N}_i$, which equals  $\sum_{j\neq i}a_{ij}$.
\item For a directed graph $\bar{\mathcal{G}}$ on the node set $\{0, 1, \cdots, N\}$, one uses a matrix $\mathcal{H}=\mathcal{L}+\mathcal{D}$ to describe it structure,
where $\mathcal{L}$ is the Laplacian matrix of its subgraph $\mathcal{G}$ and $\mathcal{D}=diag(d_1, \cdots, d_N)$ with $d_i=1$ if node $(0,i)$ is an edge of graph $\bar{\mathcal{G}}$, and with $d_i=0$ otherwise. Obviously, the structure of the graph $\bar{\mathcal{G}}$ can also be described by its Laplacian $\bar{\mathcal{L}}$.
\end{itemize}
\end{definition}

It is noted that the graph describing the interaction topology of nodes $\{0, 1, \cdots, N\}$
can vary with time. To account this we need to consider all
possible  graphs $\{\bar{\mathcal{G}}_p: p\in \mathcal {P}\}$, where $\mathcal
{P}$ is an index set for
 all  graphs  defined on  nodes $\{0,1,\cdots,N\}$. Obviously, $\mathcal{P}$ is a finite set.
We use $\{\mathcal{G}_p: p\in \mathcal {P}\}$ to denote subgraphs defined on
vertices $\{1,\cdots,N\}$. The dependence of  the graphs upon time
can be characterized by a switching law $\sigma: [0,
\infty)\rightarrow \mathcal {P}$ which is a piece-wise constant and right continuous map; that is,  at each time $t$, the
underlying graph is $\bar{\mathcal{G}}_{\sigma(t)}$.

For each agent $i\in \{1, \cdots, N\}$, if agent $j$ is a neighbor of agent
$i$,
the relative information $x_j-x_i$ is feedback to agent $i$ with a gain matrix $K$ to be design later.
The leader-following consensus problem consists of designing for each agent $i\in \{1, \cdots, N\}$ a distributed protocol which is a linear feedback, or a dynamical feedback of
\begin{eqnarray}\label{protocol}
z_i=\Sigma_{j\in \mathcal{N}_i(t)}(x_j-x_i)+d_i(t)(x_0-x_i)
 \end{eqnarray}
 such that the closed-loop systems (\ref{2.1})-(\ref{2.2}) achieve the following collected behaviors:
\begin{eqnarray}\label{leaderfollowing}
\lim_{t\rightarrow \infty}\|x_i(t)-x_0(t)\|=0, \quad i=1, \cdots, N.
\end{eqnarray}
To solve the leader-following consensus problem, the following assumption is proposed throughout this paper.
\begin{assumption}\label{stabilizable}
The pair $(A,B)$ is stabilizable.
\end{assumption}

The following result presents an averaging method for stability of fast time-varying linear systems.
The general result for nonlinear systems can be found in \cite{aeyels1999,kosut1987,bellman1985}.
For convenient of our later use, we rewrite the result in the following form.
\begin{lemma}\label{lemma2}
Consider a linear time-varying  systems $\dot x(t)=A(t)x(t)$ with $A(\cdot): \mathbb{R}\rightarrow \mathbb{R}^{n\times n}$. If  there exists an increasing sequence of times $t_k, k\in \mathbb{Z}$, with  $t_k \rightarrow +\infty$ as $k\rightarrow +\infty$,  $t_k \rightarrow -\infty$ as $k\rightarrow -\infty$, and $t_{k+1}-t_k\leq T$ for some $T>0$,  such that  $\forall t_k$, the following average systems
\begin{eqnarray}
\dot {\bar x}(t)=\bar A_k \bar x(t), \quad \bar A_k=\frac{\int_{t_k}^{t_{k+1}}A(t)dt}{t_{k+1}-t_k}, k=0,1,\cdots
\end{eqnarray}
are asymptotically stable, then there exists $\alpha^*>0$ such that the following fast time-varying system
\begin{eqnarray}
\dot x(t)=A(\alpha t)x(t)
\end{eqnarray}
is asymptotically stable for all $\alpha> \alpha^*$.
\end{lemma}
\begin{remark}\label{remark1}
It has been shown in \cite[Remark 4]{aeyels1999} that
the value $\alpha^*$  can be estimated from $T$ by solving the  equation
\begin{eqnarray}\label{alpha}
e^{\frac{KT}{\alpha}}\frac{T}{\alpha}=\frac{1}{K}\left(-1+\sqrt{1+\frac{v}{K_vKT}}\right)
\end{eqnarray}
for $\alpha$, where $T>0$ is defined above and $K_v>0, K>0, v>0$ are parameters which can be determined from the system matrix; furthermore, this equation
has for every $T>0$, $K_v>0, K>0, v>0$ exactly one positive solution $\alpha$.
Now fixing  $K_v>0, K>0, v>0$, we show that as $T\rightarrow 0$ the corresponding solution $\alpha=\alpha(T)\rightarrow 0$; indeed, $T\rightarrow 0$ rends the right hand side of  (\ref{alpha}) go to infinity, thus requiring $\frac{T}{\alpha}$ on the left hand side of  (\ref{alpha}) to go to infinity, being thus resulting in
$\alpha \rightarrow 0$.
Therefore, appropriately  choosing a small $T>0$ gives
a solution $\alpha=\alpha^*<1$.
\end{remark}

The following rank property of Kronecker product will be used. The proof is
 straightforward, being  thus  omitted.
\begin{lemma}\label{kron}
For any matrices $P, Q_1, \cdots, Q_n$ of appropriate dimensions, the following property holds:
\begin{eqnarray*}
rank({P \otimes \left(
            \begin{array}{c}
              Q_1 \\
              Q_2 \\
              \vdots\\
              Q_n
            \end{array}
          \right)})
          =rank(
\left(
  \begin{array}{c}
    P\otimes Q_1 \\
    P\otimes Q_1 \\
    \vdots\\
    P\otimes Q_n
  \end{array}
\right))
\end{eqnarray*}
\end{lemma}

The following  result will also be used later.
\begin{lemma}\label{lemma4}
Consider an $n$-order differential system $\dot x(t)=A_1x(t)+A_2y(t)$ with
$A_1\in \mathbb{R}^{n \times n},A_2\in \mathbb{R}^{n \times m}$, and $y(t)\in \mathbb{R}^m$. If $A_1$ is Hurwitz
and $\lim_{t\rightarrow \infty}y(t)=0$,  then $\lim_{t\rightarrow \infty}x(t)=0$.
\end{lemma}
{\bf Proof:}  Let $x(t,x_0,y(t))$ denote the solution of $\dot x(t)=A_1x(t)+A_2y(t)$
with initial state $x_0$ at $t=0$.
Since $A_1$ is Hurwitz, there exist positive number $\alpha,\gamma_1$ and $\gamma_2$ such that
\begin{eqnarray*}
\|x(t,x_0,y(t))\|\leq \gamma_1 \|x_0\|e^{-\alpha t}+\gamma_2 \|y(t)\|_{\infty},
\end{eqnarray*}
where $\|y(t)\|_{\infty}={\rm ess\,sup}_{t\geq 0}\|y(t)\|$.
Since $\lim_{t\rightarrow \infty}y(t)=0$, then for any $\varepsilon>0$,
there exists a $T>0$ such that
$\gamma_2\|y(t)\|<\varepsilon /2$.
Similarly, $\gamma_1 \|x_0\|e^{-\alpha t}<\varepsilon /2$. Therefore, $\|x(t,x_0,y(t))\|<\varepsilon$.
This completes the proof. \hfill $\blacksquare$


\section{Leader-Following Consensus of Multiple LTI Systems}
This section presents the leader-following consensus of multiple stablizable LTI systems under switching topology.   Unlike most results in the literature, we do not impose assumption that $A$ is neutrally stable.
For completeness, we first review a result from  \cite{ni2010} when the graph is fixed and undirected.
\begin{theorem}
For  the multi-agent system {\rm(\ref{2.1})}-{\rm(\ref{2.2})} associated with connected graph $\bar{\mathcal{G}}$ under Assumption
{\rm\ref{stabilizable}}, let  $P>0$  be a solution  to the  Riccati inequality
\begin{eqnarray}\label{riccati}
PA+A^TP-2\delta PBB^TP+I_n<0,
\end{eqnarray}
where $\delta$ is the smallest eigenvalue of the structure matrix $\mathcal H$ of graph  $\bar{\mathcal{G}}$(which is shown to be positive therein), then under the control law $u_i=Kz_i$ with $K=B^TP$
all the agents follow the leader from any initial conditions.
\end{theorem}

We now treat the leader-following consensus problem under switching topologies and directed graph case.
Denoting the state error between the agent $i$ and the leader as $\varepsilon_i=x_i-x_0$, then the dynamics of
$\varepsilon_i$ is
\begin{eqnarray*}
\dot \varepsilon_i
&=& A\varepsilon_i+Bu_i\\
&=& A\varepsilon_i+BK\sum_{j\in \mathcal{N}_i(t)}(\varepsilon_j-\varepsilon_i)-BKd_i(t)\varepsilon_i, \quad i=1,\cdots,N.
\end{eqnarray*}
By introducing
$\varepsilon=(\varepsilon_1^T, \varepsilon_2^T, \cdots,\varepsilon_N^T)^T$,
one has
\begin{eqnarray}\label{error}
\dot \varepsilon  &=&  (I_N \otimes A)\varepsilon-(I_N \otimes B) (\mathcal{L}_{\sigma(t)}\otimes I_m) (I_N \otimes K)\varepsilon-(I_N \otimes B)(\mathcal{D}_{\sigma(t)}\otimes I_m) (I_N \otimes K)\varepsilon          \nonumber\\
                  &=&  [I_N \otimes A-(\mathcal{L}_{\sigma(t)}+\mathcal{D}_{\sigma(t)})\otimes (BK)]\varepsilon   \nonumber\\
                  &=&  [I_N \otimes A-\mathcal{H}_{\sigma(t)}\otimes (BK)]\varepsilon.
\end{eqnarray}
The remaining issue is finding conditions on the switching topologies( i.e., conditions on the switching law $\sigma$) under which one can synthesize a feedback gain matrix $K$ such that the zero solution of   systems (\ref{error}) is asymptotically stable.

As treated in \cite{hong2007}, consider an infinite sequence of
nonempty,  bounded and contiguous  time intervals $[t_k, t_{k+1}), k=0,1,\cdots,$
with $t_0=0$, $t_{k+1}-t_k\leq T$ for some constant $T>0$. Suppose
that in each interval $[t_k, t_{k+1})$ there is a sequence of $m_k$
nonoverlapping subintervals
\begin{eqnarray*}
[t_k^1, t_k^2), \cdots, [t_k^j, t_k^{j+1}), \cdots, [t_k^{m_k}, t_k^{m_k+1}), \quad t_k=t_k^1, \quad t_{k+1}=t_k^{m_k+1},
\end{eqnarray*}
satisfying $t_k^{j+1}-t_k^j\geq \tau, 1\leq j\leq m_k$ for  a given constant $\tau >0$,
such that during
each of such subintervals, the interconnection topology  does not
change. That is, during each time interval $[t_k^j, t_k^{j+1})$, the
graph $\bar{\mathcal{G}}_{\sigma(t)}$ is fixed and we denote it by $\bar{\mathcal{G}}_{k_j}$.
The number $\tau$ is usually call the minimal dwell time of the graphs. The $\tau >0$
can be arbitrarily small and the existence of such an number ensures that Zero phenomena dose not happen.
During each time interval $[t_k, t_{k+1})$, some or  all  of $\bar{\mathcal{G}}_{k_j}, j=1,\cdots, m_k$ are permitted to be disconnected.
We only require the graph to be jointly connected, which is defined as follows:
\begin{definition}{\bf (Joint Connectivity)}
\begin{itemize}
\item The union of a collection of graphs is a graph whose vertex and edge sets are the unions of the vertex and edge sets of the
graphs in the collection.
\item The graphs are said to be jointly connected across the time interval $[t, t+T], T>0$ if the union of graphs
$\{\bar{\mathcal{G}}_{\sigma(s)}: s\in  [t, t+T]\}$ is connected.
\end{itemize}
\end{definition}

\begin{assumption}\label{jc}
The  graphs $\bar{\mathcal{G}}_{\sigma(t)}$ are jointly connected  across each interval $[t_k,
t_{k+1}), k=0,1,\cdots$, with their length being uniformly up-bounded by  a positive number $T$ and lower-bounded by a positive number $\tau$.
\end{assumption}

The following lemma gives a property of jointly connected graphs. When the graph is undirected, this result has been
reported in \cite{hong2007,hong2007b}. We show this result is still valid when the graph is directed; its proof is
put in the appendix.
\begin{lemma}\label{lemma1}
Let matrices
$\mathcal{H}_{1},\cdots,\mathcal{H}_{m}$ be associated with the graphs $\bar{\mathcal{G}}_{1},\cdots, \bar{\mathcal{G}}_m$ respectively. If these graphs
are jointly connected, then \\
(1) all the eigenvalues of $\sum_{i=1}^{m} \mathcal{H}_i$ have positive real parts.\\
(2) all the eigenvalues of $\sum_{i=1}^{m} \tau_i \mathcal{H}_i$ have positive real parts, where $\tau_i>0$ and $\sum_{i=1}^m \tau_i=1$.
\end{lemma}

With this, an averaging method by using Lemma \ref{lemma2} is applied to study the stability of system
(\ref{error}), whose average system during each time interval $[t_k, t_{k+1}), k=0, 1, \cdots$,  is
\begin{eqnarray}\label{averagesystem}
\dot {\bar x}=\bar A_k \bar x
\end{eqnarray} with
\begin{eqnarray*}
\bar A_k &=&\frac{\int_{t_k}^{t_{k+1}}[I_N \otimes A-\mathcal{H}_{\sigma(t)}\otimes (BK)]dt}{t_{k+1}-t_k}\\
       &=& I_N \otimes A-\bar{\mathcal{H}}_{[t_k,t_{k+1}]} \otimes (BK),
\end{eqnarray*}
where  $\bar{\mathcal{H}}_{[t_k,t_{k+1}]}=\sum_{t\in [t_k, t_{k+1})}\tau_{\sigma(t)}\mathcal{H}_{\sigma(t)}$,  $\tau_j=(t_k^{j+1}-t_k^j)/(t_{k+1}-t_k)$, $j=k_1, \cdots, k_{m_k}$.
Define by $Re\lambda_{min}(\cdot)$ the least real part of the eigenvalues of a matrix.
Define
\begin{eqnarray}\label{delta}
\bar \delta=\min \big\{\inf_{\tiny{(\tau_{k_1}, \cdots, \tau_{k_{m_k}})\in \Gamma_k}}Re\lambda_{min}
(\bar{\mathcal{H}}_{[t_k,t_{k+1})})|k=0,1,\cdots \big \},
\end{eqnarray}
where
where
\begin{eqnarray*}
\Gamma_k=\{(\tau_{k_1}, \cdots, \tau_{k_{m_k}})|\sum_{j=1}^{m_k} \tau_{k_j}, \tau \leq \tau_j < 1, j=1,\cdots, m_k\}.
\end{eqnarray*}
Noting that $Re\lambda_{min}(\bar{\mathcal{H}}_{[t_k,t_{k+1})})$ depends continuously on $\tau_{k_1}, \cdots, \tau_{k_{m_k}}$ and the set $\Gamma_k$ is compact, also by referring to  Lemma \ref{lemma2},
one has
\begin{eqnarray*}
\inf_{(\tau_{k_1}, \cdots, \tau_{k_{m_k}})\in \Gamma_k}Re\lambda_{min}
(\bar{\mathcal{H}}_{[t_k,t_{k+1})})
&=&Re\lambda_{min}(\tau_{k_1}^* \mathcal{H}_{k_1}+\cdots+\tau_{k_{m_k}}^*\mathcal{H}_{k_{m_k}} )>0,
\end{eqnarray*}
which, together with the fact that the set in (\ref{delta}) is finite due to finiteness of all graphs,   implies that
$\bar \delta$ in (\ref{delta}) is a positive number.
Then the leader-following consensus control can be achieved through the following theorem.

\begin{theorem}\label{theorem2}
For  the multi-agent system {\rm(\ref{2.1})}-{\rm(\ref{2.2})} under Assumption {\rm\ref{stabilizable}},  associated with switched graphs $\bar{\mathcal{G}}_{\sigma(t)}$  under Assumption
 {\rm\ref{jc}} with $T$ small enough, let  $P>0$ be a solution to the  Riccati inequality
\begin{eqnarray}\label{riccati2}
PA+A^TP-2 \bar \delta PBB^TP+I_n<0,
\end{eqnarray}
then under the control law $u_i=Kz_i$ with $K=B^TP$
all the agents follow the leader from any initial conditions.
\end{theorem}

{\bf Proof:} We first prove that for each $k=0,1,\cdots$, the average system
(\ref{averagesystem}) is asymptotically stable.
To this end, let $T_k \in \mathbb{R}^{N\times N}$ be
an unitary  matrix  such that $T_k \bar{\mathcal H}_{[t_k,t_{k+1}]}T^*_k=\bar \Lambda_k$ be  an upper triangular matrix with the diagonal elements $\bar\lambda_1^k, \cdots, \bar\lambda_N^k$ be the eigenvalues of matrix
$\bar{\mathcal H}_{[t_k,t_{k+1}]}$, where $T^*_k$ denote the Hermitian adjoint of matrix $T_k$. Setting  $\tilde{x}=(T_k\otimes I_n)\bar x$, (\ref{averagesystem}) becomes
\begin{eqnarray}\label{transx}
\dot{\tilde x}=(I_N \otimes A-\bar\Lambda_k \otimes BK) \tilde x.
\end{eqnarray}
The stability of (\ref{transx}) is equivalent to stability of its diagonal system
\begin{eqnarray}\label{diasystem}
\dot{\tilde x}=[I_N \otimes A-diag(\bar\lambda_1^k,\cdots, \bar\lambda_N^k) \otimes BK] \tilde x,
\end{eqnarray}
or equivalent to the stability of the following $N$ systems
\begin{eqnarray}
\dot{\tilde x}_i=(A-\bar\lambda_i^kBB^TP)\tilde x_i, \quad i=1,\cdots, N.
\end{eqnarray}
Denoting $\bar\lambda_i^k=\bar\mu_i^k+\jmath \bar \nu_i^k$, where $\jmath^2=-1$, then
\begin{eqnarray*}
&&P(A-\bar\lambda_i^kBB^P)+(A-\bar\lambda_i^kBB^TP)^*P\\
&=&P[A-(\bar\mu_i^k+\jmath \bar\nu_i^k)BB^TP]+[A-(\bar\mu_i^k+\jmath \bar\nu_i^k)BB^P]^*P\\
&=&PA+A^TP-2\bar\mu_i^k PBB^TP\\
&\leq &PA+A^TP-2 \delta PBB^TP\\
&\leq& -I<0
\end{eqnarray*}
Therefore system (\ref{averagesystem}) is globally asymptotically stable for each $k=0,1,\cdots$.

Using Lemma \ref{lemma2}, we conclude that there exists a positive $\alpha^*$ dependent of $T$,  such that $\forall \alpha > \alpha^*$, the switching system
\begin{eqnarray} \label{scale}
\dot \varepsilon(t)  = [I_N \otimes A-\mathcal{H}_{\sigma(\alpha t)}\otimes (BK)]\varepsilon(t)
\end{eqnarray}
is asymptotically stable.
According to Remark \ref{remark1}, $\alpha^*$ can be made smaller than one if we choose $T$ small enough.
Since $\alpha>\alpha^*$ is arbitrary, just pick $\alpha=1$.
That is, system (\ref{error}) is asymptotically stable, which implies that leader-following consensus is achieved.
\hfill $\blacksquare$

Although the exact value of $\bar \delta$ is hard to obtain, this difficulty can be removed as follows. Noting that for two positive parameters  $\bar{\delta}^* < \bar{\delta}$, if $P>0$  is a solution of  (\ref{riccati2})  for  parameter $\bar{\delta}^*$, then this $P$ is also a solution of  (\ref{riccati2})  for  parameter $\bar{\delta}$. Thus we can compute a positive definite matrix $P$ with a small enough parameter $\bar{\delta}^*$ which is obviously independent the global information.  This treatment has an extra advantage that it make consensus control law really distributed since the feedback gain $K=B^TP$ does not include global information.

\begin{remark}
During each interval $[t_k,t_{k+1})$, the total dwell time of the $m_k$ graphs is upper bounded by a positive number $T$,  which is required to be appropriately small to make $\alpha^*<1$. This means that the dwell time of each graph can not exceed a certain bound. However, in \cite{ni2010} the dwell time for each graph can be arbitrary since there $T$ is not constrained and can be chosen arbitrarily large.
\end{remark}

\begin{remark}
Note that in (\ref{scale}) the switching signal $\sigma(\alpha t)$ and state $\varepsilon(t)$ have different time scales, while our result is obtained for system (\ref{error}) with  $\sigma( t)$ and $\varepsilon(t)$ have the same time scale, and thus  the result in our paper is not limited to fast time switching case. This distinguishes this work from \cite{stilwell2006}.
\end{remark}

\begin{remark}
It can be seen that if  $P>0$ is a  solution to (\ref{riccati}), then $\kappa P, \kappa \geq 1$,  is also a solution to (\ref{riccati}). Indeed,
$\kappa PA+\kappa A^TP-2 \kappa ^2 \bar \delta PBB^TP+\bar \delta I_n$
$= \kappa (PA+ A^TP-2 \kappa \bar \delta PBB^TP+1/\kappa  I_n)$
$\leq  \kappa (PA+ A^TP-2  \bar \delta PBB^TP+  I_n)<0$.
Therefore, $\kappa K$ is also a stabilizing feedback matrix and the $\kappa $ can be understood as the coupling strength.
\end{remark}

\section{Observer-Based Leader-Following Consensus}
This section extends the result in last section to observer-based leader-following consensus. Consider a multi-agent system consisting of $N$ agents and a leader.
The leader agent, labeled as $i=0$,  has linear dynamics as
\begin{eqnarray}\label{leader}
\begin{array}{lllll}
\dot x_0=Ax_0,\\
y_0=Cx_0
\end{array}
\end{eqnarray}
where $y_0\in \mathbb{R}^p$ is  the output of the leader.
The dynamics of each follower agent, labeled as $i\in \{1, \cdots, N\}$,  is
\begin{eqnarray}\label{follower}
\begin{array}{llllll}
\dot x_i=Ax_i+Bu_i,\\
y_i=Cx_i
\end{array}
\end{eqnarray}
where $y_i\in \mathbb{R}^p$ is the agent $i$'s observed output information, and $u_i\in
\mathbb{R}^m$ is agent $i$'s input through which the interaction or coupling between other agents is realized.
More specifically,  $u_i$ is a  dynamical feedback of $z_i$.
In this section, we assume
 \begin{assumption}\label{sta_det}
The pair (A,B) is stabilizable, and the pair (A,C) is detectable.
\end{assumption}

The observer-based feedback controller is represented as
\begin{eqnarray}\label{observerfeedback}
\begin{array}{llllll}
{\dot{\hat{\varepsilon}}}_i=A \hat{\varepsilon}_i+K_o(\hat z_i-z_i)+Bu_i,\\
 u_i=F\hat {\varepsilon}_i,
\end{array}
\end{eqnarray}
where
\begin{eqnarray}\label{hatzi}
\hat z_i=\sum_{j\in \mathcal{N}_i}(C\hat{\varepsilon}_j-C\hat{\varepsilon}_i)+d_iC\hat{\varepsilon}_i,
\end{eqnarray}
 and the matrices $K_o$ and $K$ are
to be designed later.

\begin{remark}
The term $z_i$ in (\ref{observerfeedback}) indicates that the observer receives the output variable
information from this agent's  neighbors as input, and the term $\hat z_i$ indicates that this observer exchanges  its state with its neighboring  observers.
That is, each observer is implemented  according
to its local sensing resources. Since $z_i$ and $\hat z_i$ are local, the observer is essentially distributed,
thus then feeding  the state of each observer back to the
 corresponding agent   is again a distributed control scheme.
\end{remark}

By further introducing  the following  stacked vector
 $\hat\varepsilon=(\hat\varepsilon_1^T, \cdots, \hat\varepsilon_N^T)^T$,
  $\hat z=(\hat z_1^T, \cdots, \hat z_N^T)^T$, and by using the structure matrices of  graph $ \bar{ \mathcal{G}}_{\sigma(t)}$, one has
\begin{eqnarray}\label{hateps}
\dot {\hat\varepsilon}  &=&  (I_N \otimes A)\hat\varepsilon-[\mathcal{L}_{\sigma(t)}\otimes (K_oC)] \hat\varepsilon-[\mathcal{D}_{\sigma(t)}\otimes (K_oC)] \hat\varepsilon + \nonumber\\
&& \hspace{3cm}  [\mathcal{L}_{\sigma(t)}\otimes (K_oC)]\varepsilon+ [\mathcal{D}_{\sigma(t)}\otimes (K_oC)] \varepsilon +(I_N \otimes B)u       \nonumber\\
                      &=&  [I_N \otimes (A+BF)-\mathcal{H}_{\sigma(t)}\otimes (K_oC)]\hat\varepsilon +[\mathcal{H}_{\sigma(t)}\otimes (K_oC)][\varepsilon.
\end{eqnarray}
Then
\begin{eqnarray}\label{close}
\begin{array}{llll}
\dot { \varepsilon}=(I_N\otimes A)\varepsilon+ [I_N \otimes (BF)]\hat{\varepsilon}\\
\dot{\hat{\varepsilon}}=[\mathcal{H}_{\sigma(t)}\otimes (K_oC)]\varepsilon+[I_N\otimes A+I_N\otimes (BF)-\mathcal{H}_{\sigma(t)}\otimes (K_oC)]\hat{\varepsilon}
\end{array}
\end{eqnarray}
Let $e=\hat{\varepsilon}-\varepsilon$, that is
\begin{eqnarray*}
\left(
  \begin{array}{c}
    \varepsilon \\
    e \\
  \end{array}
\right)
=
\left(
  \begin{array}{cc}
    I_{nN} & 0 \\
    -I_{nN} & I_{nN} \\
  \end{array}
\right)
\left(
  \begin{array}{c}
    \varepsilon \\
    \hat{\varepsilon} \\
  \end{array}
\right).
\end{eqnarray*}
Under this coordinate transformation, system (\ref{close}) becomes
\begin{eqnarray}\label{system_xe}
\left(
  \begin{array}{c}
   \dot \varepsilon \\
    \dot e \\
  \end{array}
\right)
=
\left(
  \begin{array}{cc}
    I_N \otimes (A+BF) & I_N \otimes BF \\
    0 &  I_N \otimes A-\mathcal{H}_{\sigma(t)}\otimes (K_oC)\\
  \end{array}
\right)
\left(
  \begin{array}{c}
    \varepsilon \\
    e \\
  \end{array}
\right).
\end{eqnarray}

Therefore, observer-based leader-following consensus consists in, under jointly connected graph condition,  designing matrices $K_o$ and $F$ such that system (\ref{system_xe}) is asymptotically stable.
By separate principle and by referring to Lemma \ref{lemma4}, system (\ref{system_xe}) can by made asymptotically stable through carrying  out the following two-procedure design:
\begin{itemize}
\item Design matrix $K_o$ such the switched system $\dot e =[I_N \otimes A-\mathcal{H}_{\sigma(t)}\otimes (K_oC)]e$ is asymptotically stable;
\item Design matrix $F$ such the  $\dot \varepsilon =[I_N \otimes (A+BF)]\varepsilon$ is asymptotically stable;
\end{itemize}
The first step can be realized by referring Theorem \ref{theorem2}, replacing the pair $(A, B)$ with $(A^T,C^T)$.
The second step is a standard state feedback control problem.
We  summarize above analysis  in the following theorem. The rest of its proof is essentially similar
to that of  Theorem \ref{theorem2}, and is thus omitted for save of space.
\begin{theorem}\label{theorem4}
Consider the multi-agent systems (\ref{leader}-\ref{follower}) associated with  switching graphs $\bar{\mathcal{G}}_{\sigma(t)}$ under the Assumptions \ref{jc},
\ref{sta_det} with $T$ small enough.  Let  $P>0$ be a solution  to the  Riccati inequality
\begin{eqnarray}\label{riccati}
PA^T+AP-2\bar \delta PC^TCP+I_n<0,
\end{eqnarray}
then under the control law {\rm(\ref{observerfeedback})} with $K_o=PC^T$ and $F$ being such that $A+BF$ is Hurwitz,
all the agents follow the leader from any initial conditions.
\end{theorem}

\section{Simulation Results}
\begin{figure}[h]
\centering
\includegraphics[width=11cm]{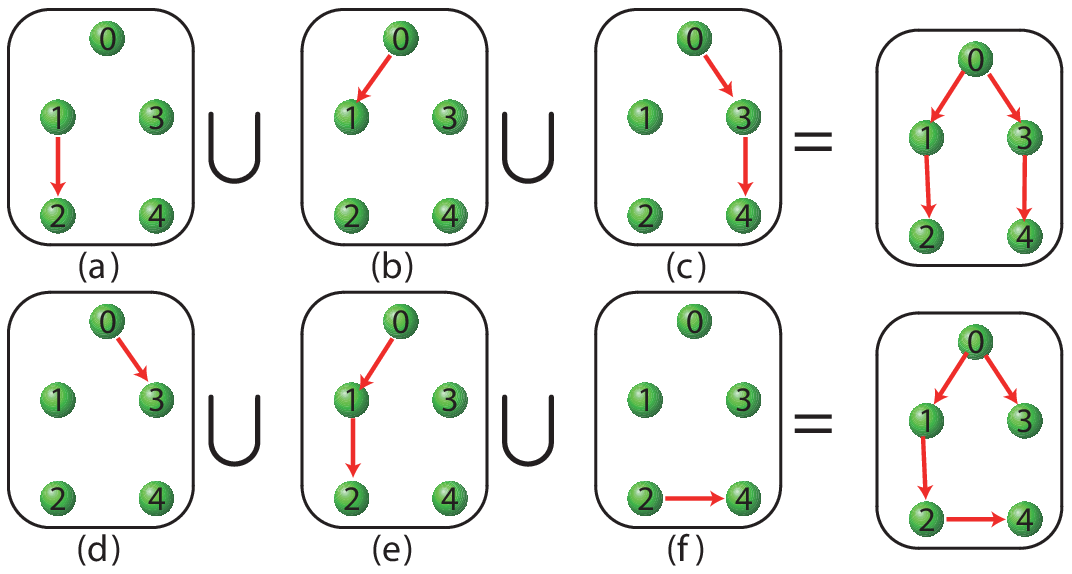}
\caption{Six possible interaction topologies between the leader and the agents.}
\label{topology}
\end{figure}

\begin{figure}[h]
\centering
\includegraphics[width=11cm]{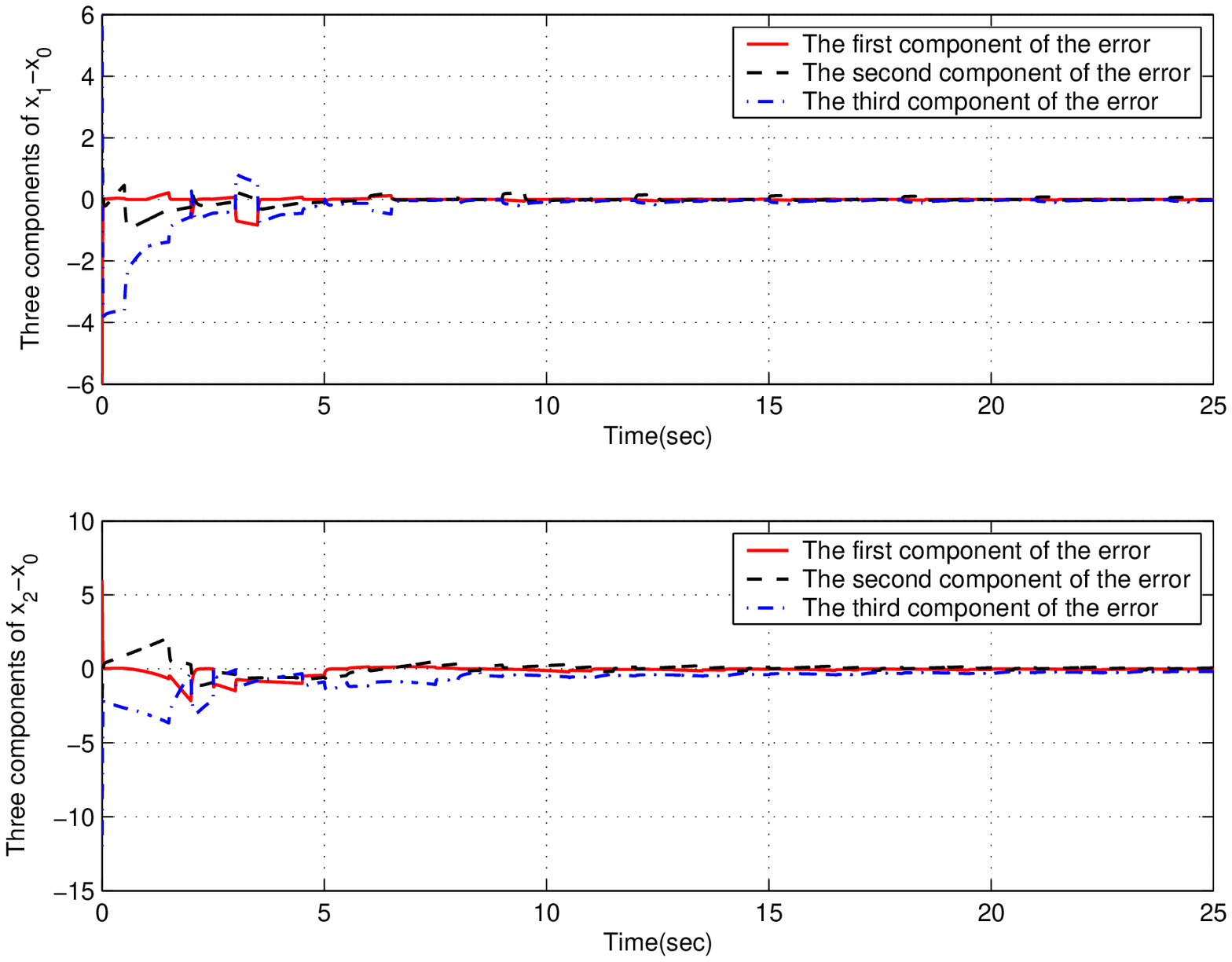}
\includegraphics[width=11cm]{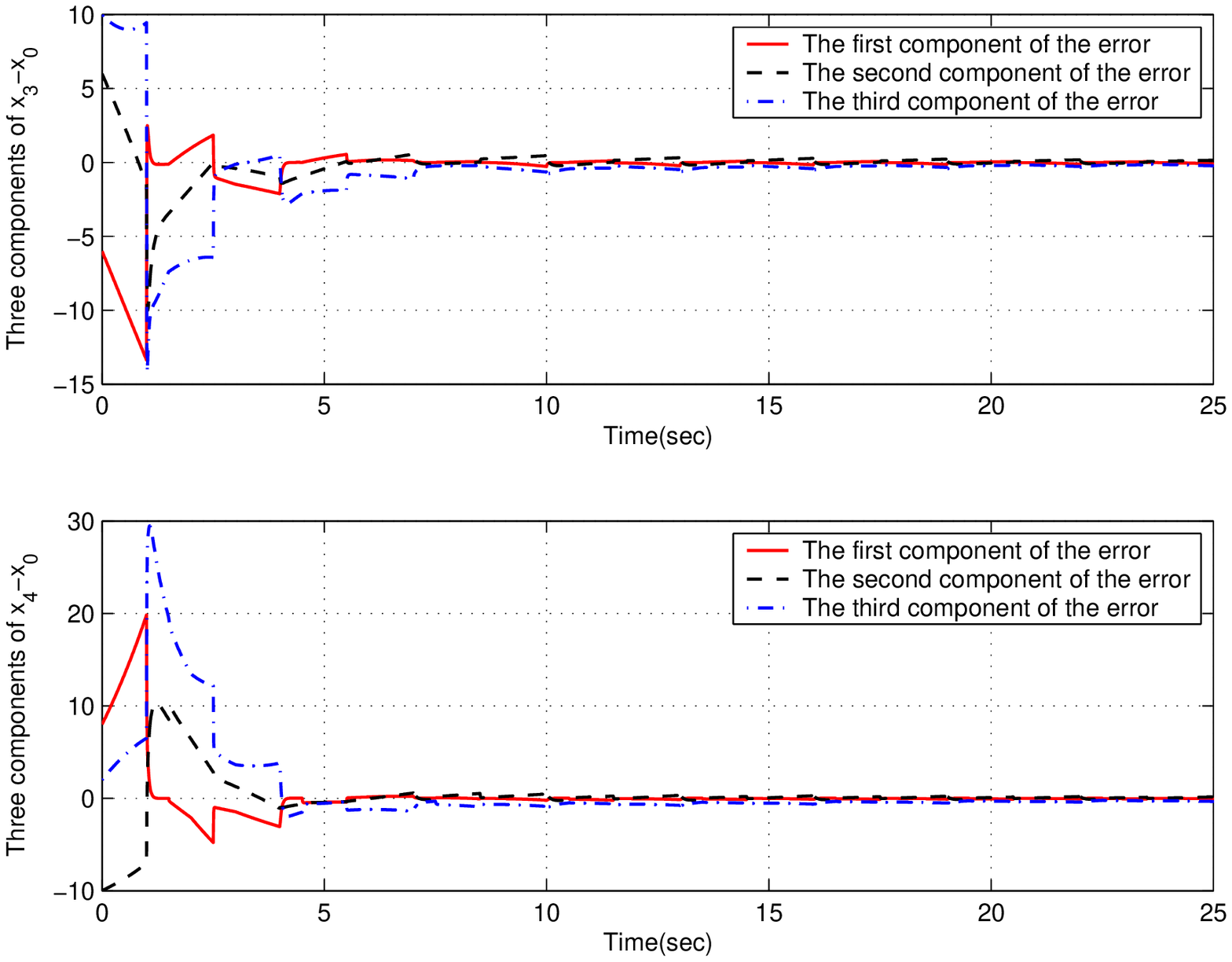}
\caption{Simulation for Theorem \ref{theorem2}: The error trajectories between the leader and each agent.}
\label{figth2}
\end{figure}

\begin{figure}[h]
\centering
\includegraphics[width=11cm]{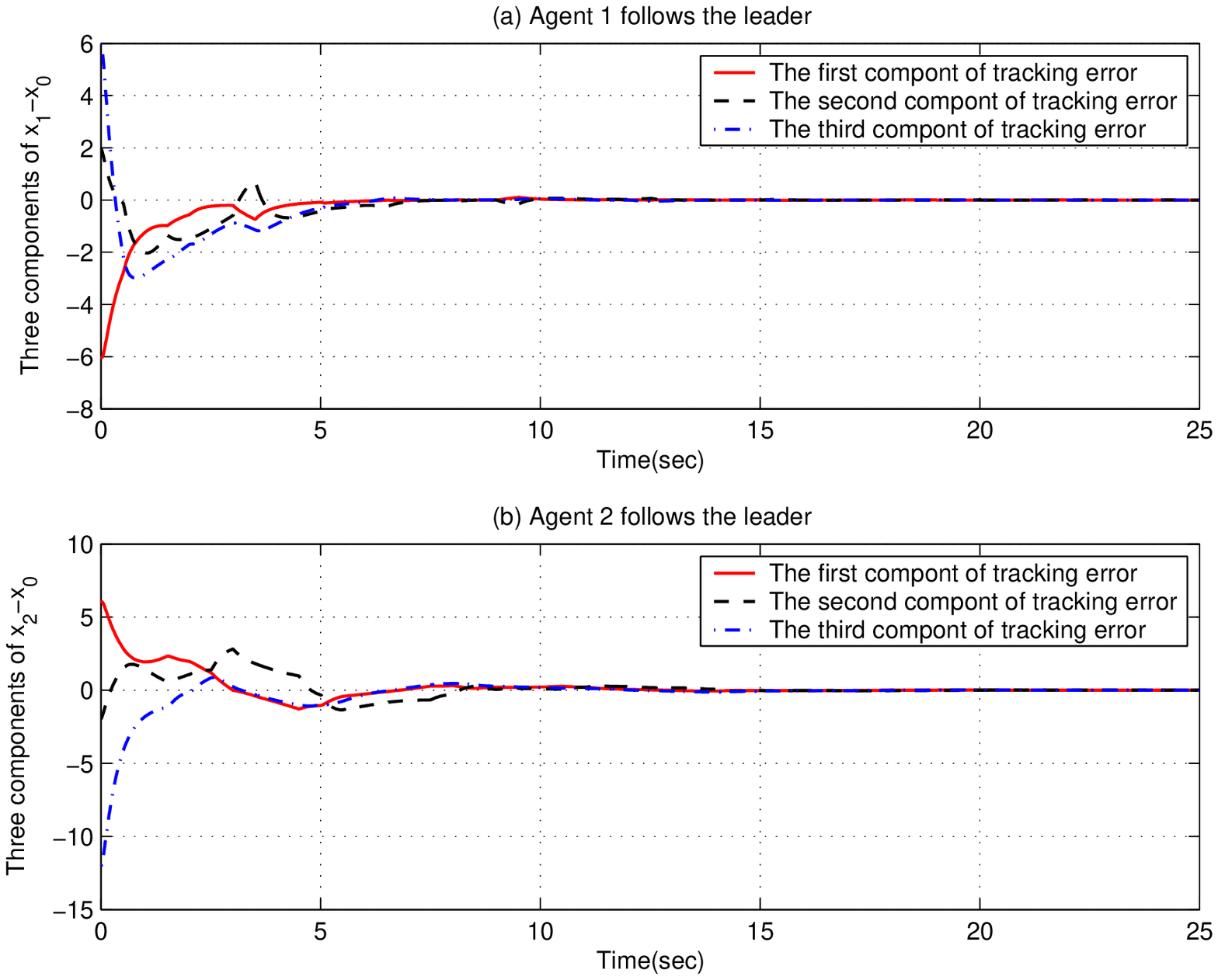}
\includegraphics[width=11cm]{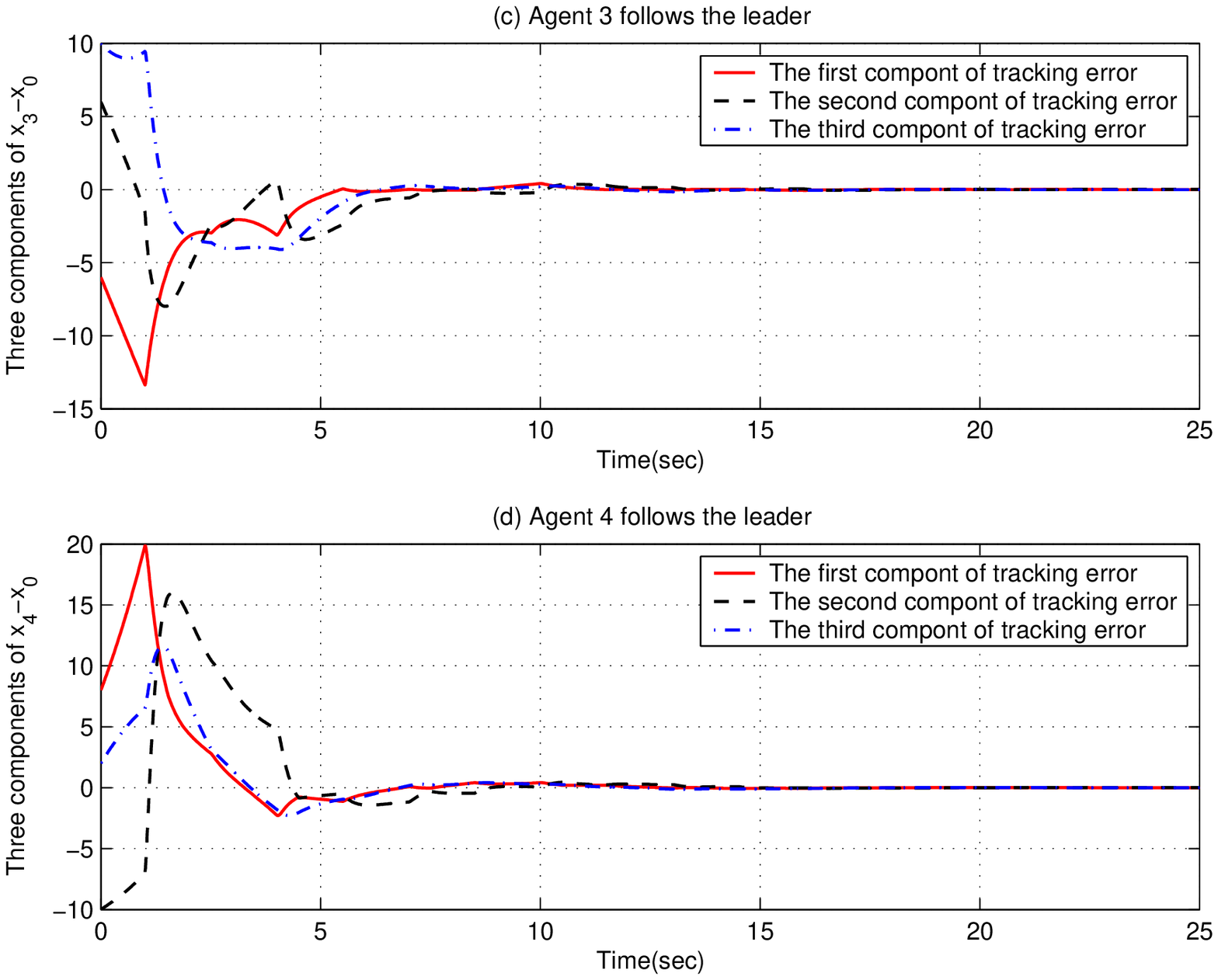}
\caption{Simulation for Theorem \ref{theorem4}: The error trajectories between the leader and each agent.}
\label{figth4}
\end{figure}

In this section, we give two examples to illustrate the validity
of the results. Consider a multi-agent system consisting of a leader and four
agents. Assume the system matrices are
\begin{eqnarray*}
A=\left(
    \begin{array}{ccc}
    0.5548 &  -0.5397 &  -0.0757\\
    0.3279 &  -0.0678 &  -0.4495\\
   -0.0956 &  -0.6640 &   0.0130
    \end{array}
  \right),
B=\left(
    \begin{array}{cc}
    3  &   5\\
     3 &   -2\\
    -8 &   -8
    \end{array}
  \right),
C=\left(
    \begin{array}{ccc}
      1  &  -1  &   2\\
    -4   &  2   & -3
    \end{array}
  \right)
\end{eqnarray*}
We suppose that  possible
interaction graphs are
$\{\bar{G}_1,\bar{G}_2,\bar{G}_3,\bar{G}_4,\bar{G}_5,\bar{G}_6\}$
which are shown in Figure {\rm\ref{topology}}, and the
interaction graphs are switching as
$\bar{G}_1\rightarrow\bar{G}_2\rightarrow\bar{G}_3\rightarrow\bar{G}_4
\rightarrow\bar{G}_5\rightarrow\bar{G}_6\rightarrow
\bar{G}_1\rightarrow \cdots $, and each graph is active for $1/2$
second. Since the graphs $\bar{G_1}\cup \bar{G_2}\cup
\bar{G_3}$ and $\bar{G_4}\cup \bar{G_5}\cup \bar{G_6}$ are
connected, we can choose $t_k=k, t_{k+1}=k+3/2$ and $t_k^0=k,
t_k^1=k+1/2, t_k^2=k+2/2, t_k^3=k+3/2$ with $k=0,1\cdots$.
We choose a small parameter
$\bar \delta=1/3 min(0.3820, 0.1732)=0.0057$.  The matrices $K$ in Theorem \ref{theorem2} and $K_O, F$ Theorem \ref{theorem4} are calculated as
\begin{eqnarray*}
K=\left(
    \begin{array}{ccc}
     0.7520 &   5.9852 &  -2.7041\\
   12.6966  & -3.8441  &  1.6419
    \end{array}
  \right)
\end{eqnarray*}
and
\begin{eqnarray*}
F=\left(
    \begin{array}{ccc}
    0.6338 &  -0.5087 &   0.3731\\
   -0.9077 &   0.4509  & -0.1938
    \end{array}
  \right),
K_O=\left(
    \begin{array}{ccc}
    -6.7092 &   -9.1532\\
   -9.6111  &  4.1353\\
    7.7514  &  -1.1756

    \end{array}
  \right),
\end{eqnarray*}
respectively.
With the same initial condition, the simulation results of Theorem 2 and Theorem 4 are shown in Figure \ref{figth2}
and Figure \ref{figth4}, respectively.

\section{Conclusions}

This paper presents an averaging   approach to leader-following consensus problem of multi-agent with linear dynamics, without imposing neutrally stable condition on agent dynamics.
The interaction topology is switching.
The proposed protocols force the follower agents to follow
the independent leader trajectory.
The result is extended to observer-based protocols design.
Such design can be separated as a two-step procedure: design asymptotically stable distributed
 observers and asymptotically stable observer-state-feedback protocols.

\section*{Appendix}
The appendix is devoted to the proof of Lemma \ref{lemma1}. To this end, we first cite the  following result.

\begin{lemma}(Ger\v{s}gorin)\cite{horn1985}
For any matrix $G=[g_{ij}]\in \mathbb{R}^{N\times N}$, all the eigenvalues of $G$ are located in the union of $N$ Ger\v{s}gorin
discs
\begin{eqnarray*}
Ger(G):=\cup_{i=1}^{N}\{z\in C: |z-g_{ii}|\leq \sum_{j\neq i}|g_{ij}|\}.
\end{eqnarray*}
\end{lemma}

The definition of weighted graph will also be used in the proof of what follows.
If we assign each edge $(i,j)$ of graph $\bar{\mathcal{G}}$ a weight $w_{ij}$,  we obtain a weighted graph $\bar{\mathcal{G}}_{\mathcal W}=(\bar{\mathcal{V}}, \bar{\mathcal{E}},\bar{\mathcal W})$, where $\bar{\mathcal{W}}=[w_{ij}]$.
For an graph $\bar{\mathcal{G}}$ and any positive number $k>0$, the graph $k\bar{\mathcal{G}}$ is defined to be a weighted graph obtained from $\bar{\mathcal{G}}$ by assigning a weight $k$ to each existing edge of $\bar{\mathcal{G}}$.
For two graphs $\bar{\mathcal{G}}_1$ and $\bar{\mathcal{G}}_2$, their  union is a weighted graph  and the weight
 for edge $(i,j)$  is the sum of weights for  the two edges $(i,j)$ in the graphs $\bar{\mathcal{G}}_1$ and $\bar{\mathcal{G}}_2$ respectively.
The weighted Laplacian  of graph $\bar{\mathcal{G}}_{_{\mathcal{W}}}$
is  defined as $\bar{\mathcal{L}}_{_{\mathcal{W}}}=-\bar{\mathcal{A}}_{_{\mathcal{W}}}+\bar{\Lambda_{_{\mathcal{W}}}}$, where $\bar{\mathcal{A}}_{_{\mathcal{W}}}=[w_{ij}a_{ij}]$ and $\bar{\Lambda_{_{\mathcal{W}}}}(i,i)=\sum_{j\neq}w_{ij}a_{ij}$;
the weighted structure matrix of graph $\bar{\mathcal{G}}_{_{\mathcal{W}}}$
is  defined as $\mathcal{H}_{_{\mathcal{W}}}=\mathcal{L}_{_{\mathcal{W}}}+\mathcal{D}_{_{\mathcal{W}}}$, where
$\mathcal{L}_{_{\mathcal{W}}}$ is the weighted Laplacian of the subgraph $\mathcal{G}_{_{\mathcal{W}}}$ of $\bar{\mathcal{G}}_{_{\mathcal{W}}}$, and
$\mathcal{D}_{_{\mathcal{W}}}=diag(w_{_{01}}d_1,\cdots, w_{_{0N}}d_N)$.

{\bf Proof of Lemma \ref{lemma1}:}  (1) Denote the Laplacian matrix  and the structure matrix of graph $\bar{\mathcal{G}}_i$ by  $\bar{\mathcal{L}}_i$ and $\mathcal{H}_i$ respectively,  and  denote the Laplacian matrix   of graph $\mathcal{G}_i$ by $\mathcal{L}_i$.

We first prove the case when $m=1$. By definitions, it can be easily verified  the following relationship
\begin{eqnarray*}
\bar{\mathcal{L}}_1=\begin{pmat}({|..})
0 & 0 & \cdots & 0 \cr\-
-d_1 &  &  &  \cr
\vdots &  & \mathcal{H}_1&  \cr
-d_N &  &  &   \cr
\end{pmat}.
\end{eqnarray*}
Since the graph $\bar{\mathcal{G}}_1$ is connected, then $rank(\bar{\mathcal{L}}_1)=N$ \cite{ren2005}. Thus the sub-matrix $\bar{\mathcal {M}}_1$  formed by the last $N$ rows of $\bar{\mathcal{L}}_1$ has rank $N$.
Note that
\begin{eqnarray*}
\left(
  \begin{array}{c}
    -d_1 \\
    \vdots \\
    -d_N \\
  \end{array}
\right)
=
\mathcal {D}_1
\left(
  \begin{array}{c}
    1 \\
    \vdots \\
    1 \\
  \end{array}
\right)
=
\mathcal {L}_1\left(
  \begin{array}{c}
    1 \\
    \vdots \\
    1 \\
  \end{array}
\right)
+
\mathcal {D}_1
\left(
  \begin{array}{c}
    1 \\
    \vdots \\
    1 \\
  \end{array}
\right)
=
\mathcal {H}_1
\left(
  \begin{array}{c}
    1 \\
    \vdots \\
    1 \\
  \end{array}
\right),
\end{eqnarray*}
that is, the first column of matrix $\bar{\mathcal {M}}_1$ is a linear combination of its last $N$ columns.
Therefore, $rank(\mathcal {H}_1)=N$.   Furthermore, we claim the eigenvalues of  the matrix  $\mathcal {H}_1$ are located in closed-right half plan; indeed,  from Ger\v{s}gorin theorem, all the eigenvalues of $H$ are located in
\begin{eqnarray*}
Ger(H)=\cup_{i=1}^{N}\left\{z\in \mathbb{C}: |z-l_{ii}-d_i|\leq |N_i|\right\},
\end{eqnarray*}
and therefore they are located in the closed-right half plan by noting that $l_{ii}=|N_i|$ and $d_i\geq 0$.
We thus conclude that all the eigenvalues of $\mathcal {H}_1$ have positive real parts.

We proceed to prove the case when $m>1$. Obviously,  for a union $\bar{\mathcal{G}}_{_U}$ of a group of weighted graphs $\{\bar{\mathcal{G}}_1, \cdots, \bar{\mathcal{G}}_m\}$,
its weighted Laplacian matrix  $\bar{\mathcal{L}}_{_U}$ is the sum of the  Laplacian matrices  $\{\bar{\mathcal{L}}_1, \cdots, \bar{\mathcal{L}}_m\}$ of graphs $\{\bar{\mathcal{G}}_1, \cdots, \bar{\mathcal{G}}_m\}$, and
\begin{eqnarray*}
\bar{\mathcal{L}}_U=\begin{pmat}({|..})
0 & 0 & \cdots & 0 \cr\-
-d_1^U &  &  &  \cr
\vdots & \mathcal{H}_1+ & \cdots & +\mathcal{H}_m  \cr
-d_N^U &  &  &   \cr
\end{pmat},
\end{eqnarray*}
where
\begin{eqnarray*}
\left(
  \begin{array}{c}
    -d_1^U \\
    \vdots \\
    -d_N^U \\
  \end{array}
\right)=
\left(
  \begin{array}{c}
    -d_1^1 \\
    \vdots \\
    -d_N^1 \\
  \end{array}
\right)+\cdots+
\left(
  \begin{array}{c}
    -d_1^m \\
    \vdots \\
    -d_N^m \\
  \end{array}
\right)
\end{eqnarray*}
with $(d_1^j, d_2^j, \cdots, d_N^j)^T$ be the diagonal elements of matrix $\mathcal {D}_j, j=1,\cdots, m$.
When the graphs are jointly connected, that is, when the $\bar{\mathcal{G}}_{_U}$ is connected, the matrix
$\bar{\mathcal{L}}_{_U}$ has a simple zero eigenvalue. Argue in a manner similar to that of $m=1$ case, it can be shown that the all the eigenvalues of the matrix $\mathcal{H}_1+  \cdots  +\mathcal{H}_m$ have positive real parts.

(2) Similar discussion as given in (1) for the weighted graphs $\tau_1 \bar{\mathcal{G}}_1, \cdots, \tau_N \bar{\mathcal{G}}_N$ yields the conclusion.

\section*{ACKNOWLEDGEMENT}
\small
This work is  supported by  the NNSF of China(61104096, 11101203, 60904024), the Youth Foundation of Jiangxi Provincial Education Department of China(GJJ12132), JXNSF(20114BAB201002) and SDNSF(ZR2011FQ014).

\makesubmdate

%
%
%
%
%
%
%

\makecontacts

\end{document}